\documentclass[a4paper,11pt,twoside]{article}
\usepackage[a4paper,top=5cm,bottom=5cm,left=4cm,right=4cm]{geometry}

\usepackage[latin1]{inputenc}
\usepackage[T1]{fontenc}
\usepackage[english]{babel}

\usepackage{sansmath}
\usepackage{libertine}
\usepackage{libertinust1math}
\usepackage[scaled]{beramono}

\usepackage{calrsfs}
\usepackage{enumitem}

\setlist{nolistsep}
\usepackage{mathtools}
\usepackage{mathrsfs}
\usepackage{mfirstuc}
\usepackage{multicol}

\usepackage[overload]{textcase}
\usepackage{tasks}
\settasks{counter-format=tsk[1].}
\usepackage{bbm}
\usepackage{amsmath}
\usepackage{amssymb}
\usepackage{amsfonts}
\usepackage{amsthm}
\usepackage{wasysym}
\usepackage{dsfont}
\usepackage{graphicx}
\usepackage{wrapfig}
\usepackage[outdir=./]{epstopdf}
\usepackage[labelfont={bf,it},textfont=it]{caption}
\usepackage[hidelinks]{hyperref}
\usepackage{widetable}
\usepackage[scr=pxtx, scrscaled=1]{mathalfa}
\usepackage{xfrac}
\usepackage{soul}
\usepackage{tikz}
\usepackage{setspace}
\usepackage{yfonts}
\usepackage{url}

\theoremstyle{plain}
\theoremstyle{definition}
\newtheorem{definition}{\sffamily Definition}[section]

\newtheorem{lemma}[definition]{\sffamily Lemma}

\newtheorem{proposition}[definition]{\sffamily Proposition}

\newtheorem{theorem}{\sffamily Theorem}
\theoremstyle{remark}
\newtheorem*{remark}{\normalfont\sffamily Remark}

\newtheorem*{example*}{\normalfont\sffamily Example}

\usepackage{sectsty}
\allsectionsfont{\sffamily\large}
\subsectionfont{\sffamily}
\newcommand*\myrule[1][.25\textwidth]{%
	\tikz {\path [fill, draw] (0,0) [out=0, in=180] to +(.5*#1,1pt) [out=0, in=180] to +(.5*#1,-1pt) [out=180, in=0] to +(-.5*#1,-1pt) [out=180, in=0] to cycle;}}

\makeatletter
\def\@maketitle{%
  \newpage
  \null
  \vskip 2em%
  \begin{center}%
  \let \footnote \thanks
    \begin{minipage}{.7\textwidth}
    	\begin{center}
    		\sffamily\Large\textbf{{ \@title }}\par
    	\end{center}
    \end{minipage}
    \vskip 1.5em%
    {\large
      \lineskip .5em%
      \vskip -.8em
      \vspace{-1em} 
      \begin{tabular}[t]{c}%
        {\@author}
      \end{tabular}\par}%
    \vskip .5em%
    \par
  \vskip .7em
  \vspace{-1.5em} 
  \myrule[.15\textwidth]
  \vskip 1.5em
  \end{center}%
  }
\makeatother

\renewenvironment{abstract}{%
\vskip -2em
\hfill\begin{center}\small
\begin{minipage}{0.75\textwidth}}
{\par\noindent
\end{minipage}
\end{center}
\vskip 2.5em}
\makeatletter

\usepackage{caption}
\usepackage{authblk}

\usepackage{fancyhdr}
\fancyheadoffset{0cm}

\usepackage{blindtext}
\usepackage{parskip}

\newcommand{\blu}[1]{#1}

\newcommand{\blusec}[1]{#1}

\newcommand{\ensemblenombre}[1]{\mathbb{#1}}

\newcommand{\Z}{\ensemblenombre{Z}}

\newcommand{\R}{\ensemblenombre{R}}
\newcommand{\Dim}{d}
\newcommand*{\1}{\textrm{{\usefont{U}{fplmbb}{m}{n}1}}}
\newcommand{\A}{\mathcal{A}_R}
\newcommand{\Ao}{\mathcal{A}^a}
\newcommand{\HardOpen}{U} 
\newcommand{\Intensity}{z}
\newcommand{\Parameters}{\Intensity, \beta}

\newcommand{\Conf}{\omega}
\newcommand{\BoundaryConf}{\gamma}
\newcommand{\ConfSpace}{\Omega}
\newcommand{\SigmaAlgebra}{\mathcal{F}}
\newcommand{\Poisson}{\pi}

\newcommand{\Hamiltonian}{H}
\newcommand{\Potential}{\phi}
\newcommand{\Specification}{\mathscr{P}}
\newcommand{\PartitionFunction}{Z}

\newcommand{\Gibbs}{\mathcal{G}_{\Parameters}(\Potential)}
\newcommand{\GibbsDis}{\mathcal{G}^{a\Z^d}_{\Parameters}(\Potential)}
\newcommand*{\HCdiam}{R}
\newcommand{\Width}{a}
\newcommand{\SuperLambda}{\Lambda}
\newcommand{\abs}[1]{\lvert #1\rvert}
\newcommand{\Borel}{\mathcal{B}(\R^\Dim)}
\newcommand*{\CRuelle}{C_\phi(\beta)}

\usepackage{authblk}

\makeatletter
\def\namedlabel#1#2{\begingroup
    #2%
    \def\@currentlabel{#2}%
    \phantomsection\label{#1}\endgroup
}
\makeatother

\newcommand{\defeq}{%
	\mathrel{ \vcenter{\baselineskip0.5ex \lineskiplimit0pt \hbox{\scriptsize.} \hbox{\scriptsize.}} } =}

\begin{document}

\title{An explicit Dobrushin uniqueness region for Gibbs point processes with repulsive interactions}

\author[1]{Pierre Houdebert}
\author[2]{Alexander Zass}
\affil[1]{Universit\"at Potsdam}
\affil[2]{Universit\"at Potsdam, WIAS Berlin}


\maketitle

\begin{abstract}
We present a uniqueness result for Gibbs point processes with interactions that come from a non-negative pair potential; in particular, we provide an explicit uniqueness region in terms of activity $z$ and inverse temperature $\beta$. The technique used relies on applying to the continuous setting the classical Dobrushin criterion. We also present a comparison to the two other uniqueness methods of cluster expansion and disagreement percolation, which can also be applied for this type of interactions.
\bigbreak

\medbreak
\end{abstract}


\section{Introduction}
\label{sec:introduction}
In this work we present a uniqueness result for a class of Gibbs point processes in $\R^\Dim$, $d\geq 2$, where the Gibbsian interaction is given by a non-negative pair potential.

In statistical mechanics, Gibbs point processes are often defined through the \emph{Do\-bru\-shin--Lanford--Ruelle (DLR) equations}
\cite{Dobrushin_1968, Lanford_Ruelle_1969}
that prescribe their conditional laws. This large class of point processes is widely used, as it allows for various types of interaction: among others, it can be in the form of a $k$-body potential; depend on geometric features like the Delaunay tessellation or the Area interaction processes; be attractive or repulsive (for an introduction to Gibbs point processes see \cite{Georgii1979, dereudre_mini_cours}).

Setting aside the -- far from trivial -- existence problem, which is not the focus of this paper, the natural question is that of uniqueness (or lack thereof) of Gibbs point processes. Indeed, in the setting of Gibbs point processes that are characterised by a pair potential $\Potential$, an activity parameter $\Intensity > 0$ and an inverse temperature $\beta \geq 0$, we place particular focus on obtaining an explicit uniqueness region of the parameters $\Intensity,\beta$.
Heuristically, it is expected that for fixed $\beta$, uniqueness is achieved for activities $\Intensity$ small enough.
We remark however that such a behaviour has actually been disproved for the specific case of the Widom--Rowlinson model with random radii with heavy tails, see \cite{Dereudre_Houdebert_2019_JSP_PhaseTransitionWR}.

Already for lattice systems, the uniqueness question is one of main interest in the community, and different arguments and methods exist, including Peierls' argument \cite{peierls_1936}, the Dobrushin criterion from \cite{Dobrushin_1968}, cluster expansion (see \cite{mayer_1937,ruelle_livre_1969}), a characterisation due to Lebowitz and Martin L\"of \cite{lebowitz_loef_1972}, and disagreement percolation (see, for example, \cite{georgii_haggstrom_maes, vandenberg_1993, vandenberg_maes}).
No one method is a priori stronger than the others: it is a question of finding which is better adapted to the specific model one wishes to prove uniqueness of the Gibbs measure for.

This consideration also holds true in the continuous framework, where the uniqueness question is even more delicate. The three techniques we analyse in the setting of this paper -- namely  the \emph{Dobrushin contraction criterion}, \emph{cluster expansion}, and \emph{disagreement percolation} -- all work under different assumptions, and yield different parameter domains in which uniqueness holds. It is therefore generally complicated to compare their efficacy.

In the main result of this paper, Theorem \ref{thm:uniqueness_ours}, we provide a set of (simple to test) assumptions which lead to a uniqueness result for small activity. The strength of this result lies in the explicit nature of the uniqueness domain it yields, i.e. the following parameter region for $\Parameters$:
\begin{equation*}
	\Big\{(\Parameters)\in\R_{>0}\times\R_{\geq 0}\colon \Intensity < \Big(\sup_{x \in \R^\Dim} \ \int_{\R^\Dim} e^{-\beta \Potential(x,y)}dy\Big)^{-1}\Big\},
\end{equation*}
and in the simplicity of its proof, which makes use of the Dobrushin criterion. The celebrated original uniqueness criterion from \cite{Dobrushin_1968} -- very general, but presented only for discrete Gibbs models -- can be applied to continuum models by first decomposing the space $\R^\Dim$ into disjoint cubes of some side length $\Width>0$. Letting this discretisation parameter $\Width$ tend to $0$ \blu{(in \eqref{eq:finalbound})}, we are able to obtain the explicit uniqueness bound above.

\blusec{We prove our result for potentials with a hard-core component close to the origin (Assumption (A1) in Subsection \ref{sec:Dobrushin_thm_ours}), as such a requirement allows us to only consider configurations having at most one point in each small cube, which in turn simplifies the computations when taking the limit as the size of the cubes goes to $0$.
However, since the probability of having more than one point in a small cube vanishes together with its size, we conjecture that Theorem \ref{thm:uniqueness_ours} could still be valid without this hard-core assumption. See Subsection \ref{sec:discussion_assumptions} for more detailed comments.}

The second aspect of our work is a comparison between the three uniqueness techniques; in particular, in Subsection \ref{sec:comparison}, we compare it to the explicit uniqueness regions that can be obtained from existing works on cluster expansion and disagreement percolation, in \cite{jansen_2019} and \cite{Hofer-temmel_Houdebert_2018}, respectively. What transpires from this comparison is that Theorem \ref{thm:uniqueness_ours} yields a larger uniqueness region than what can be obtained via cluster expansion \cite{Fernandez_Procacci_Scoppola_2007,jansen_2019}.
Furthermore, as expected, for $\beta$ small enough, the result is also better than the one obtained from disagreement percolation, yielding a larger range of possible activities $z$ for which uniqueness holds.

\medbreak
The article is organised as follows:  
in Section \ref{sec:preliminaries} we introduce the formalism used in this work.
In Subsection \ref{sec:Dobrushin_thm_ours} we introduce the assumptions needed and state the uniqueness theorem, which is then proved in Subsection \ref{sec:proof}.
In Subsection \ref{sec:discussion_assumptions} we comment on the assumptions and give possible generalisations to our work.
In Subsections \ref{sec:compareDob} and \ref{sec:comparison} we discuss the optimality of our result, and compare it to existing results coming from cluster expansion and disagreement percolation.

\section{The setting} \label{sec:preliminaries}
\subsection{Configuration space}
In this work we consider point configurations in $\R^\Dim$, $\Dim\geq 2$.

More precisely, we endow $\R^\Dim$ with the usual Euclidean distance $\abs{\cdot}$ and Borel $\sigma$-algebra $\Borel$, and set the configuration space $\ConfSpace$ to be the set of \emph{locally finite configurations} $\Conf$ on $\R ^\Dim$, i.e. measures of the form $\omega = \sum_i \delta_{x_i}$, with $\# \Conf_{\Lambda}< \infty$ for any bounded $\Lambda$ in $\Borel$ (here $\#$ denotes the cardinality of the support of the configuration, and $\Conf_{\Lambda}$ the restriction of $\omega$ to $\Omega_\Lambda$).

As the configurations we consider here are \emph{simple}, i.e. with no overlapping points, we also denote a configuration $\omega = \sum_i \delta_{x_i}$ by the subset of $\R^\Dim$ on which it is supported: $\omega = \{x_1,\dots,x_n,\dots\}$. Consequently, $\omega_\Lambda = \omega\cap\Lambda$. We write $\Conf' \Conf \defeq \Conf' \cup \Conf$ for the \emph{concatenation} (or union) of two configurations.

We endow $\ConfSpace$ with the usual $\sigma$-algebra $\SigmaAlgebra$ generated by the counting functions on bounded Borel sets, $\omega\mapsto \#\omega_\Lambda$.

For any $\Lambda \subset \R^\Dim$, $\ConfSpace_{\Lambda}\subset \ConfSpace$ denotes the subset of configurations supported on $\Lambda$ (and by $\SigmaAlgebra_{\Lambda}$ the corresponding $\sigma$-algebra).

Let $\abs{\cdot}$ be the usual $\Dim$-dimensional Lebesgue measure. We denote by $v_\Dim\defeq \abs{B(0,1)}$ the volume of the unit ball in $\R^\Dim$.
\medbreak
On the space $\ConfSpace$, we consider the probability measure $\Poisson^{\Intensity}$ given by the distribution of the homogeneous \emph{Poisson point process} with intensity $\Intensity >0$. Recall that this means that:
\begin{enumerate}[label=\roman*.]
\item For every bounded set $\Lambda\subset\R^\Dim$, the distribution of the number of points in $\Lambda$ under $\Poisson^\Intensity$ is a Poisson distribution of mean $\Intensity \abs{\Lambda}$;
\item Given the number of points in a bounded set $\Lambda$, the said points are independent and uniformly distributed in $\Lambda$.
\end{enumerate}
For any bounded set $\Lambda \subset \R^\Dim$, we denote by $\Poisson^{\Intensity}_\Lambda$ the restriction of $\Poisson^{\Intensity}$ to $\ConfSpace_\Lambda$.

For more details on Poisson point processes see, for example, \cite{daley_vere_jones,last_penrose_2017}.

\subsection{Interactions and Gibbs point processes}
\label{sec:interaction}
As usual, we add an interaction on the Poisson point process by considering the notion of Gibbs specifications associated to a given Hamiltonian. More precisely, let $\Potential$ be a symmetric \emph{non-negative} (i.e. repulsive) \emph{pair potential}
\begin{equation*}
\Potential : \
\R^\Dim \times \R^\Dim \to \R_{\geq 0} \cup \{+\infty\},
\end{equation*}
and define, for any bounded set $\Lambda \subset \R^\Dim$, the \emph{$\Lambda$-Hamiltonian} by setting
\begin{equation*}
\Hamiltonian_\Lambda (\Conf)
\defeq
\displaystyle 
\sum_{\substack{\{x,y\} \subset \Conf 
\\ \{x,y\}\cap \Lambda \not = \emptyset } }
\Potential (x,y),\quad \Conf\in\ConfSpace.
\end{equation*}
Since the potential $\Potential$ is non-negative,
this quantity is well-defined for any configuration $\Conf\in\ConfSpace$. In Section \ref{sec:Dobrushin_thm_ours} we provide the more precise assumptions that are needed for the main result of this work, but we remark now that we do not assume translation invariance.

\begin{definition}
	The \emph{Gibbs specification} associated to the potential $\Potential$ on a bounded measurable set $\Lambda \subset \R^\Dim$, with boundary condition $\BoundaryConf$, activity $z>0$ and inverse temperature $\beta\geq 0$, is given by the following probability measure on the configurations with support in $\Lambda$:
\begin{equation*}
    \Specification_{\Lambda, \BoundaryConf}^{\Parameters}(d \Conf_\Lambda) \defeq \frac{ e^{-\beta \Hamiltonian_{\Lambda} (\Conf_{\Lambda} \BoundaryConf_{\Lambda^c}) }} {\PartitionFunction^{\Parameters}_{\Lambda,\BoundaryConf}}  \ \Poisson^{\Intensity}_{\Lambda}(d\Conf_{\Lambda}),\quad \Conf_\Lambda\in\ConfSpace_\Lambda,
\end{equation*}
where the normalisation factor is the \emph{partition function}
\begin{equation*}
    \PartitionFunction^{\Parameters}_{\Lambda,\BoundaryConf} \defeq \int_{\ConfSpace} e^{-\beta \Hamiltonian_{\Lambda} (\Conf_{\Lambda} \BoundaryConf_{\Lambda^c})}\Poisson^{\Intensity}_{\Lambda}(d\Conf_{\Lambda}) \in (0,1].
\end{equation*}
\end{definition}
Since $\PartitionFunction^{\Parameters}_{\Lambda,\BoundaryConf}$ is finite, the Gibbs specification is always \emph{well-defined}.

\begin{definition}\label{def_gibbs_measures}
A probability measure $P$ on $\ConfSpace$ is said to be a \emph{Gibbs point process} associated to the potential $\phi$, with activity $\Intensity>0$ and inverse temperature $\beta\geq 0$, denoted $P \in \Gibbs$, if for every bounded measurable function $f$ and for all bounded Borel $\Lambda \subset \R^\Dim$, the following hold
\begin{equation}
\label{eq:dlr_continuous}\tag{DLR}
    \int_{\ConfSpace} f \ d P = \int_{\ConfSpace} \int_{\ConfSpace} f(\Conf_{\Lambda}  \BoundaryConf_{\Lambda^c} )\ \Specification^{\Parameters}_{\Lambda, \BoundaryConf}(d \Conf_\Lambda)\ P (d \BoundaryConf), 
\end{equation}
called \emph{DLR equations} after Dobrushin, Lanford and Ruelle, and which prescribe the conditional probabilities of a Gibbs point process.
\end{definition}

The first question that arises in Gibbs point processes theory is whether there exists at least one solution to the DLR equations.
This important and difficult problem has been studied for many different interactions and settings: from the the classical works of \cite{ruelle_livre_1969} and \cite{pechersky_zhukov_1999}, to more recent works considering the case of geometrical interactions (\cite{dereudre_2009,dereudre_drouilhet_georgii}), infinite range pair potentials (\cite{dereudre_vasseur_2020}), and unbounded interactions in the context of marked point processes (\cite{roelly_zass_2020}).

The existence of such a measure in the setting of this paper is a known result (for example by D. Ruelle in \cite{ruelle_1970}), and is therefore not the subject of our work.
For completeness, we state it here:
\begin{proposition}
Let $\Potential$ be a non-negative and symmetric pair potential. Let $\beta\geq 0$ and assume that
\begin{equation*}
    \forall x\in\R^\Dim,\quad \int_{\R^\Dim}  \left( 1 - e^{-\beta \Potential (x,y)} \right) dy <\infty.
\end{equation*}
Then, for any activity $z>0$, there exists at least one Gibbs point process $P\in\Gibbs$.
\end{proposition}

In what follows, we first prove a simple and explicit uniqueness criterion derived from the standard Dobrushin technique, and then compare it to criteria coming from the two other techniques of cluster expansion and disagreement percolation.

We do not explore the topic of phase transition; we only mention that the question of non-uniqueness of the Gibbs point process is of major interest and very few results are known. In particular, the existing literature mainly deals with coloured (multi-species) models like the \emph{Widom-Rowlinson} model; see for instance \cite{chayes_kotecky, Dereudre_Houdebert_2019_JSP_PhaseTransitionWR}.
In these works, phase transition is proved by showing that one species ``dominates'' the others when the activity of the points is large enough.

\section{Uniqueness of the Gibbs point process}
In Subsections \ref{sec:Dobrushin_thm_ours} and \ref{sec:proof}, we derive a simple and explicit uniqueness region, by applying the discrete Dobrushin contraction criterion from \cite{Dobrushin_1968} through a \emph{discretisation} parameter $\Width$, and then considering the limit as $\Width$ goes to $0$.

The remainder of the paper deals with the natural questions that arise from this result. In particular, in Subsection \ref{sec:discussion_assumptions} we discuss the assumptions of Theorem \ref{thm:uniqueness_ours}, as well as possible generalisations that could be the subject of future works. In Subsection \ref{sec:compareDob} we perform a numerical study to show that taking the limit for the mesh size $a\rightarrow 0$ yields a larger uniqueness region than that of the Dobrushin criterion for any fixed $a>0$.
Finally, in Subsection \ref{sec:comparison} we compare this uniqueness region to the ones coming from other approaches, namely \emph{Cluster expansion} and \emph{Disagreement percolation}.

\subsection{Dobrushin uniqueness region}
\label{sec:Dobrushin_thm_ours}

We present here the assumptions that are required in the statement of our main result. See Subsection \ref{sec:discussion_assumptions} for some comments on these conditions.

\noindent \textbf{Assumptions.}
In what follows we assume that the (not necessarily translation-invariant) pair potential $\Potential$ is such that the following three conditions are satisfied:
\begin{enumerate}
\item[(A1)] 
\blu{\emph{Repulsive and hard-core close to the origin:} $\Potential$ is non-negative, and there exists a uniform hard-core diameter $R>0$ such that}
\begin{equation*}
    \forall\, x,y\in\R^\Dim\colon \abs{x-y}< R,\quad \Potential(x,y)=+ \infty.
\end{equation*}
\item[(A2)]
\emph{Uniform regularity of the potential:}
\begin{equation*}
	\CRuelle\defeq \sup_{x \in \R^\Dim} \int_{\R^\Dim} 
\left( 1 - e^{-\beta \Potential (x,y)} \right) dy 
<+\infty.
\end{equation*}
\end{enumerate}
For the third assumption we first need to introduce some notations.
Let $\Width >0$, and partition the space $\R^\Dim$ into cubes of side length $a$, centred in the points of the lattice: for any $i\in \Width \Z^\Dim$, these are
\begin{equation}\label{eq:lambdai}
	\Lambda_{\Width,i} \defeq i + \left(-\frac{\Width}{2} , \frac{\Width}{2}\right]^\Dim,
\end{equation}
then define the following ``local supremum'' of the Mayer function
\begin{equation}\label{eq:Psi}
	\Psi_\Width(x,y)\defeq \sup_{\bar{y} \in \Lambda_{\Width,i} } \left( 1 - e^{- \beta \Potential (x,\bar{y})}\right) \quad \text{if } y\in\Lambda_{\Width,i}.
\end{equation}
The last assumption is given by the following
\begin{enumerate}
\item[(A3)] \emph{Regularity of the induced Mayer function:}
\begin{equation*}
\begin{split}
	\sup_{x \in \R^\Dim} \int_{\R^\Dim} \left( 1 - e^{-\beta \Potential (x,y)} \right) dy &= \lim_{\Width \to 0} \ \sup_{x \in \R^\Dim} \ \int_{\R^\Dim} \Psi_\Width(x,y) \ dy\\
	&\blusec{= \inf_{0<\Width<2R/\sqrt{\Dim}} \ \sup_{x \in \R^\Dim} \ \int_{\R^\Dim} \Psi_\Width(x,y) \ dy}.
\end{split}
\end{equation*}
\end{enumerate}

We are now ready to state the main result of this paper:
\begin{theorem}
\label{thm:uniqueness_ours}
Let the pair potential $\Potential$ satisfy assumptions (A1), (A2), (A3). Furthermore, assume
\begin{equation}\label{eq:th_uniqueness_region}
	\Intensity < \CRuelle^{-1}.
\end{equation}
Then there exists a unique Gibbs point process $P\in\Gibbs$.
\end{theorem}

\subsection{Proof of Theorem \ref{thm:uniqueness_ours}}
\label{sec:proof}
In this subsection we provide the proof to the above theorem. As already stated in the introduction, the result is an application to the continuous setting of the classical Dobrushin technique for lattice models. To do this, we use a standard discretisation technique, via a parameter $a>0$ that defines the mesh size $a$. The novelty of our result, and what leads to the explicit uniqueness region, is to consider the Dobrushin criterion in the limit as $a\to 0$; this final step of the proof is done in Subsection \ref{subsec:bound}.

\subsubsection{The Dobrushin contraction method in the lattice}
The proof of Theorem \ref{thm:uniqueness_ours} relies on the classical Dobrushin criterion \cite{Dobrushin_1968}.
In this subsection we describe the setting and the results as they apply to our model (for a general presentation see, for example, \cite{georgii_livre}).

Let $\Xi$ be a complete separable metric space, which we call the \emph{spin space}. Fix $\Width>0$, and consider the \emph{discrete configuration} space $\Xi^{\Width \Z^\Dim}$, equipped with the standard cylinder $\sigma$-algebra. Let $\Pi = \left(\Pi_\Gamma(\cdot\vert\xi)\right)_{\Gamma \subset \Width \Z^\Dim, \xi \in  \Xi^{ \Width \Z^\Dim}}$ be a \emph{discrete specification}, i.e. a consistent family of conditional probability measures indexed by a finite $\Gamma \subset \Width \Z^\Dim$ and a discrete configuration $\xi \in  \Xi^{ \Width \Z^\Dim}$. Furthermore, for any event $A$, $\Pi_\Gamma(A\vert\xi)$ only depends on the restriction $\xi_{\Gamma^c}$ of $\xi$ to $\Xi^{\Gamma^c}$.
\begin{definition}
	A probability measure $Q$ on $\Xi^{\Width\Z^\Dim}$ is said to be a \emph{Gibbs measure compatible with the specification $\Pi$} if, for any finite $\Gamma \subset \Width \Z^\Dim$ and any bounded measurable function $g$, it satisfies
	\begin{equation*}
		\int g(\xi) Q(d\xi) = \int g(\xi'_\Gamma \ \xi_{\Gamma^c} ) \Pi_\Gamma(d \xi'_\Gamma\vert\xi) Q(d \xi).
	\end{equation*}
\end{definition}

The original result, proved by Dobrushin in \cite{Dobrushin_1968} reads
\begin{theorem}
\label{thm:unique_dobru}
Let $\Width>0$ and $\Pi_i \defeq \Pi_{\{i \}}$. If
\begin{equation*}
    \sup_{i \in \Width \Z^\Dim} \ \sum_{j\not = i} \sup_{\substack{\xi,\tilde\xi\in\Xi^{\Width\Z^\Dim}\\\xi_k = \tilde{\xi}_k\ \forall k\neq j}} d_{TV}\left( \Pi_i\ (\cdot\vert\xi),\Pi_i\ (\cdot\vert\tilde{\xi}) \right) < 1,
\end{equation*}
then there exists at most one Gibbs measure on $\Xi^{\Width\Z^\Dim}$ compatible with the specification $\Pi$.
\end{theorem}

An easy generalisation, which we use in the following subsection to restrict the set of allowed configurations, is the following
\begin{lemma}\label{lemma:dobrushin}
	Let $\Width>0$ and consider $A\subset\Xi^{a\Z^\Dim}$ such that $Q(A)=1$ for any Gibbs measure $Q$ on $\Xi^{\Width\Z^\Dim}$ compatible with the specification $\Pi$. If
\begin{equation}\label{eq:dobrushin_crit}
    \sup_{i \in \Width \Z^\Dim} \ \sum_{j\not = i} \sup_{\substack{\xi,\tilde\xi\in A\\\xi_k = \tilde{\xi}_k\ \forall k\neq j}} d_{TV}\left( \Pi_i\, (\cdot\vert\xi),\Pi_i\, (\cdot\vert\tilde{\xi}) \right) < 1,
\end{equation}
then there exists at most one Gibbs measure in $\Xi^{a\Z^\Dim}$ compatible with the specification $\Pi$.	
\end{lemma}

\begin{remark}
	As seen in the above Lemma, the hard-core assumption (A1) allows us to restrict the set of possible boundary conditions in \eqref{eq:dobrushin_crit}. This is also achieved by Do\-brushin and Pecherski in \cite{dobrushin_pecherski_1983} -- extended to the continuous setting in \cite{pechersky_zhukov_1999,belitsky_pechersky_2002} -- but while the criterion provides a uniqueness result for more general interactions than those considered here, it does not generally provide an explicit uniqueness region, which is one of the goals of our work. Indeed, as seen in Lemma 6 of \cite{pechersky_zhukov_1999}, their method yields uniqueness for small (non-explicit) values of the parameters.
\end{remark}

\subsubsection{Correspondence between continuous and lattice models}
In order to apply Lemma \ref{lemma:dobrushin}, one must express the continuous model as a lattice model. The representation we make use of here does not lose any of the information from the continuous model, so that the uniqueness properties of the two models are indeed equivalent (see \cite{dobrushin_pecherski_1983,pechersky_zhukov_1999}). In this it differs from the cell-gas model presented by A. Rebenko in \cite{rebenko_cellgas}, where the idea is that uniqueness for the discrete model implies uniqueness for the continuous one, but it is not a one-to-one correspondence.

Fix $\Width>0$, and consider again the partition of $\R^\Dim$ into the cubes $\Lambda_{\Width, i}$ defined in \eqref{eq:lambdai}. Set the spin space to be the space of configurations supported on the closure of $\Lambda_{\Width,0}$, $\Xi \defeq \ConfSpace_{\bar\Lambda_{\Width,0}}$, endowed with its $\sigma$-algebra $\SigmaAlgebra_{\bar\Lambda_{\Width,0}}$. The lattice configurations are then subsets $\xi\subset\Xi^{\Width\Z^\Dim}$. The correspondence between the two models is given by the following measurable embedding
\begin{equation*}
\begin{split}
    T\colon &\ConfSpace \to \Xi^{\Width \Z^\Dim}\\
    &\Conf \mapsto \xi = (\xi_j)_{j\in\Width\Z^\Dim},
\end{split}
\end{equation*}
defined by setting, for all $j\in \Width\Z^\Dim$, 
\begin{equation*}
	\xi_j \defeq \Conf_{\Lambda_{\Width, j}} - j
= \{ x- j: x \in \Conf_{\Lambda_{\Width, j}} \}\subset \Lambda_{a,0}.
\end{equation*}
\blu{Its inverse $T^{-1}$ can be naturally extended from $T(\ConfSpace)$ to the whole of $\Xi^{\Width\Z^\Dim}$ by considering 
\begin{equation*}
	T^{-1}(\xi) = \bigcup_{j\in\Width\Z^\Dim} \xi_j +j.		
\end{equation*}
\begin{remark}
	A lattice configuration $\xi\in\Xi^{\Width\Z^\Dim}\setminus T(\ConfSpace)$ contains points in $\bar \Lambda_{a,0}\setminus \Lambda_{a,0}$ and so $T^{-1}$ maps it to a continuous configuration $\gamma$ with overlapping points. Nevertheless, such configurations have probability $0$ under the Gibbs point process as $\phi(x,x) = +\infty$.
\end{remark}}
We next define the lattice specification by setting, for any finite $\Gamma \subset \Width \Z^\Dim$,
\begin{equation*}
	\Pi_{\Gamma}(\cdot\vert\xi) \defeq
\Specification^{\Parameters}_{\Lambda,\BoundaryConf},
\end{equation*}
where $\Lambda = \cup_{i \in \Gamma} \Lambda_{\Width, i}$
and $\BoundaryConf = T^{-1}(\xi)$, and let $\GibbsDis$ be the set of probability measures on $\Xi^{\Width \Z^\Dim}$ compatible with the specification $\Pi$. Thanks to the hard-core assumption (A1) and the previous remark, it can be seen, as in \cite{pechersky_zhukov_1999}, that for $\Width< 2R/\sqrt{\Dim}$, the map $T$ induces a bijective map
\begin{equation*}
    \hat{T}: \Gibbs \to \GibbsDis.
\end{equation*}
\blu{The Gibbs point processes in $\Gibbs$ are supported on the set of \emph{$R$-hard-core configurations} (commonly, \emph{allowed}) $\A\subset\ConfSpace$, i.e. $\BoundaryConf\in\A$ is such that if $\{x,y\}\subset\gamma$, then $\abs{x-y}\geq R$; the Gibbs measures in $\GibbsDis$ are supported on their image via $T$, i.e. $A\defeq T(\A)$.}

We can then restrict our study to only the allowed configurations $\A$. Indeed, thanks to Lemma \ref{lemma:dobrushin} and the one-to-one correspondence between continuous and lattice models, in order to prove uniqueness of the original Gibbs point process, it is enough to show that, for some $a>0$,
\begin{equation}\label{eq:dobrushin_cont}
	\sup_{i\in \Width \Z ^\Dim} \sum_{j\in \Width \Z ^\Dim\setminus\{i\}} k^{(a)}_{i,j} < 1,
\end{equation}
where we have set, for any $\{i,j\} \subset \Width \Z^\Dim$,
\begin{equation*}
    k^{(a)}_{i,j} \defeq \sup_{\substack{\BoundaryConf,\tilde\BoundaryConf\in\A\\\BoundaryConf_{\Lambda_{\Width,j}^c} = \tilde{\BoundaryConf}_{\Lambda_{\Width,j}^c}}} d_{TV} \left( \Specification^{\Parameters}_{\Lambda_{\Width, i}, \BoundaryConf}; \Specification^{\Parameters}_{\Lambda_{\Width, i}, \tilde{\BoundaryConf}} \right).
\end{equation*}

\subsubsection{Computation of the coefficients}

We consider $0<\Width<2R/\sqrt{d}$ so that the cube of side length $a$ is contained in the hard core: $\Lambda_{\Width,0} \subset B(0,\HCdiam)$. Thanks to assumption (A1), this implies that every allowed configuration 
$\BoundaryConf\in\A$ has at most one point in each cube $\Lambda_{\Width,k},\ k\in\Width\Z^\Dim$.

Fix $i\in \Width\Z^\Dim$, and for the sake of simplicity, let $\SuperLambda$ denote the cube $\Lambda_{a,i}$. For any $j \in \Width \Z^\Dim\setminus\{i\}$, consider two configurations that differ only inside the cube $\Lambda_{\Width,j}$, i.e.
\begin{equation}\label{eq:confSup}
	\BoundaryConf, \tilde{\BoundaryConf}\in\A \colon \BoundaryConf_{\Lambda^c_{\Width,j}} =\tilde{\BoundaryConf}_{\Lambda^c_{\Width,j}}.
\end{equation}
Without loss of generality, we can assume that
${\PartitionFunction^{\Parameters}_{\SuperLambda}(\tilde{\BoundaryConf})}
\leq
{\PartitionFunction^{\Parameters}_{\SuperLambda}({\BoundaryConf)}}$. Computing the total variation distance (see, e.g. Lemma 8.2.1 in \cite{resnick_2014}), we have
\begin{equation*}
\begin{split}
	& d_{TV} \left(\Specification^{\Parameters}_{\SuperLambda, \BoundaryConf}, \Specification^{\Parameters}_{\SuperLambda, \tilde{\BoundaryConf}} \right) = \sup_{B\subset\A}\abs{\Specification^{\Parameters}_{\SuperLambda, \BoundaryConf}(B) - \Specification^{\Parameters}_{\SuperLambda, \tilde{\BoundaryConf}}(B)}\\
    &= \int_{\A} \left[ \frac{e^{- \beta H_{\SuperLambda}( \Conf_{\SuperLambda} \BoundaryConf_{\SuperLambda^c})}} {\PartitionFunction^{\Parameters}_{\SuperLambda}(\BoundaryConf)} - \frac{e^{- \beta H_{\SuperLambda}( \Conf_{\SuperLambda} \tilde{\BoundaryConf}_{\SuperLambda^c})}} {\PartitionFunction^{\Parameters}_{\SuperLambda}(\tilde{\BoundaryConf})} \right]^+ \Poisson^{\Intensity}_{\SuperLambda}(d\Conf_{\SuperLambda}),
\end{split}	
\end{equation*}
where $y^+ \defeq \max(y,0)$ denotes the \emph{positive part} of $y\in\R$.
Since the configuration $\Conf_\SuperLambda$ constains at most one point in $\SuperLambda$, we obtain
\begin{equation*}
    d_{TV} \left( \Specification^{\Parameters}_{\SuperLambda, \BoundaryConf} , \Specification^{\Parameters}_{\SuperLambda, \tilde{\BoundaryConf}} \right) = \Intensity e^{-\Intensity \abs{\SuperLambda}} \int_{\SuperLambda} \left[ \frac{e^{- \beta H_{\SuperLambda}( \{x\} \cup \BoundaryConf_{\SuperLambda^c}) } } {\PartitionFunction^{\Parameters}_{\SuperLambda}(\BoundaryConf)} - \frac{e^{- \beta H_{\SuperLambda}( \{x\} \cup \tilde{\BoundaryConf}_{\SuperLambda^c}) }  } {\PartitionFunction^{\Parameters}_{\SuperLambda}(\tilde{\BoundaryConf})} \right]^+ dx.
\end{equation*}
\blu{Since $e^{-\Intensity \abs{\SuperLambda}} \leq \PartitionFunction^{\Parameters}_{\SuperLambda}(\tilde{\BoundaryConf}) \leq \PartitionFunction^{\Parameters}_{\SuperLambda}(\BoundaryConf) \leq 1$, we have that}
\begin{equation*}
\begin{split}
	e^{-\Intensity \abs{\SuperLambda}} &\left[ \frac{e^{- \beta H_{\SuperLambda}( \{x\} \cup \BoundaryConf_{\SuperLambda^c}) } } {\PartitionFunction^{\Parameters}_{\SuperLambda}(\BoundaryConf)} - \frac{e^{- \beta H_{\SuperLambda}( \{x\} \cup \tilde{\BoundaryConf}_{\SuperLambda^c}) }  } {\PartitionFunction^{\Parameters}_{\SuperLambda}(\tilde{\BoundaryConf})} \right]^+\\
	 &\leq e^{-\Intensity \abs{\SuperLambda}} \left[ \frac{e^{- \beta H_{\SuperLambda}( \{x\} \cup \BoundaryConf_{\SuperLambda^c}) } } {\PartitionFunction^{\Parameters}_{\SuperLambda}(\BoundaryConf)} - \frac{e^{- \beta H_{\SuperLambda}( \{x\} \cup \tilde{\BoundaryConf}_{\SuperLambda^c}) }  } {\PartitionFunction^{\Parameters}_{\SuperLambda}({\BoundaryConf})} \right]^+\\
	& \leq \left[ e^{- \beta H_{\SuperLambda}( \{x\} \cup \BoundaryConf_{\SuperLambda^c}) } - e^{- \beta H_{\SuperLambda}( \{x\} \cup \tilde{\BoundaryConf}_{\SuperLambda^c})} \right]^+
\end{split}
\end{equation*}
and we can now estimate
\begin{equation}\label{eq:computation_dtv_removing_boundary}
    d_{TV} \left( \Specification^{\Parameters}_{\SuperLambda, \BoundaryConf} , \Specification^{\Parameters}_{\SuperLambda, \tilde{\BoundaryConf}} \right) \leq \Intensity  \int_{\SuperLambda} \left[ e^{- \beta H_{\SuperLambda}( \{x\} \cup \BoundaryConf_{\SuperLambda^c}) } - e^{- \beta H_{\SuperLambda}( \{x\} \cup \tilde{\BoundaryConf}_{\SuperLambda^c})} \right]^+ dx.
\end{equation}
\blu{Since $\phi$ is non-negative, the right-hand side is largest when $\BoundaryConf = \emptyset$; hence from \eqref{eq:confSup}, $\tilde{\BoundaryConf} = \tilde{\BoundaryConf}_{\Lambda_{a,j}}$, with at most one point in the cube $\Lambda_{a,j}$. We then have that the following uniform bound holds over all pairs of allowed boundary configurations:}
\begin{equation}\label{eq:tvestimate}
    \forall \BoundaryConf, \tilde{\BoundaryConf}\in\A: \BoundaryConf_{\SuperLambda^c} = \tilde{\BoundaryConf}_{\SuperLambda^c}, \ d_{TV} \left( \Specification^{\Parameters}_{\SuperLambda, \BoundaryConf} , \Specification^{\Parameters}_{\SuperLambda, \tilde{\BoundaryConf}} \right) \leq \Intensity  \sup_{y \in \SuperLambda } \int_{\SuperLambda} \left( 1 - e^{- \beta \Potential (x,y) } \right) dx.
\end{equation}
In particular, taking the supremum over the configurations $\BoundaryConf, \tilde{\BoundaryConf}\in\A$ such that $\BoundaryConf_{\Lambda_{a,j}^c} = \tilde{\BoundaryConf}_{\Lambda_{a,j}^c}$ yields
\begin{equation*}
	k^{(a)}_{i,j}\leq \Intensity \sup_{y \in \Lambda_{a,j} }\int_{\Lambda_{a,i}}\left(1 - e^{- \beta \Potential (x,y) } \right) dx.
\end{equation*}

\subsubsection{Obtaining the bound}\label{subsec:bound}
Summing the coefficients $k^{(a)}_{i,j}$ over all $j\in\Width\Z^\Dim\setminus\{i\}$, leads to
\begin{equation*}
\begin{split}	
    \sum_{\substack{j \in \Width \Z^\Dim \\ j\not = i}} k^{(a)}_{i,j} & \leq \Intensity \sum_{\substack{j \in \Width \Z^\Dim \\ j\not = i}} \sup_{y \in \Lambda_{a,j} } \ \int_{\Lambda_{\Width, i}} \left( 1 - e^{- \beta \Potential (x,y) } \right) dx \\ & \leq \Intensity \int_{\Lambda_{\Width, i}} \sum_{\substack{j \in \Width \Z^\Dim \\ j\not = i}} \sup_{y \in \Lambda_{a,j} } \ \left( 1 - e^{- \beta \Potential (x,y) } \right) dx \\
    & = \Intensity \int_{\Lambda_{\Width, i}} \sum_{\substack{j \in \Width \Z^\Dim \\ j\not = i}} \frac{1}{\abs{\Lambda_{\Width, j}}} \int_{\Lambda_{\Width, j}} \sup_{\bar{y} \in \Lambda_{a,j} } \ \left( 1 - e^{- \beta \Potential (x,\bar{y}) } \right) dy\, dx \\
    & =  \frac{\Intensity}{\abs{\Lambda_{\Width, i}}} \int_{\Lambda_{\Width, i}} \int_{\R^\Dim \setminus \Lambda_{\Width, i}} \Psi_\Width(x,y)\, dy\, dx \\
    & \leq \Intensity \sup_{x \in \Lambda_{\Width, i}} \int_{\R^\Dim } \Psi_\Width(x,y)\, dy,
\end{split}
\end{equation*}
where $\Psi_\Width$ is the function defined in \eqref{eq:Psi}.
Taking the supremum over all $i \in \Width \Z^\Dim$ yields
\begin{equation}\label{eq:finalbound}
\begin{split}	
    &\sup_{i \in \Width \Z^\Dim} \sum_{\substack{j \in \Width \Z^\Dim \\ j\not = i}} \sup_{\BoundaryConf, \tilde{\BoundaryConf}\in\A } \ d_{TV} \left( \Specification^{\Parameters}_{\Lambda_{\Width, i}, \BoundaryConf} , \Specification^{\Parameters}_{\Lambda_{\Width, i}, \tilde{\BoundaryConf}} \right)\\
    & \leq \Intensity \sup_{x \in \R^\Dim} \int_{\R^\Dim } \Psi_\Width(x,y)\, dy \ \xrightarrow{\Width \to 0} \Intensity \sup_{x \in \R^\Dim} \int_{\R^\Dim } \left( 1 - e^{-\beta \Potential (x,y)} \right) dy,
\end{split}
\end{equation}
where the last convergence follows from Assumption (A3).
This means, in particular, that if $\Intensity \sup_{x \in \R^\Dim} \int_{\R^\Dim } \left( 1 - e^{-\beta \Potential (x,y)} \right) dy<1$, there exists $a>0$ such that the Dobrushin condition \eqref{eq:dobrushin_cont} is satisfied: $\sup_{i\in \Width \Z ^\Dim} \sum_{j\in \Width \Z ^\Dim\setminus\{i\}} k^{(a)}_{i,j} < 1$. This concludes the proof of Theorem \ref{thm:uniqueness_ours}.

\begin{remark}
	The uniqueness region obtained for fixed $a>0$, i.e.
	\begin{equation*}
		\Intensity \sup_{x \in \R^\Dim} \int_{\R^\Dim } \Psi_\Width(x,y)\, dy<1	
	\end{equation*}
	is not explicit due to the nature of the function $\Psi_a$. For the same reason, comparing it for different values of $a$ is not an easy feat. See the comparison in Subsection \ref{sec:compareDob}.
\end{remark}

\subsection{Discussion about the assumptions and possible generalisations}\label{sec:discussion_assumptions}
The proof of Theorem \ref{thm:uniqueness_ours} relies on the assumptions on the interaction made above, namely that it is coming from a non-negative pair potential $\Potential$ which is hard-core close to the origin (A1), with an integrability assumption (A2), and a technical regularity assumption (A3). We comment here on these requirements.

\subsubsection*{\blusec{Discussion about the assumptions}}

\blu{Firstly, we restrict our study to \emph{non-negative} potentials $\Potential$. We used this assumption to simplify the estimation of the total-variation distance in \eqref{eq:tvestimate}, since in this case $e^{-\beta H_\Lambda(\gamma)}$ is (uniformly) bounded from above by $1$. Extending our result to more general pair potentials, e.g. stable and regular like those considered in the Kirkwood--Salsburg approach of \cite{ruelle_livre_1969}, remains an open question.
In the setting of repulsive potentials we can, however, easily apply and compare all three uniqueness methods.}

Secondly, the regularity assumption (A3) is purely technical, resulting from taking $\Width \to 0$ in the proof.
The first equality is satisfied whenever the set of discontinuity of $\Potential(x,.)$ is of measure $0$, which is the case, for example, for the radial potentials with hard core presented in Subsection \ref{sec:comparison}. The result of \blusec{Theorem \ref{thm:uniqueness_ours} is obtained without making use of the second equality, which ensures instead that the uniqueness region is optimal within the Dobrushin criterion and it is consistent with numerical results (see Subsection \ref{sec:compareDob}).}

Thirdly, the regularity assumption (A2) is quite standard, and seems unavoidable when using the Dobrushin criterion.
Furthermore, a similar assumption is required when using the cluster expansion technique.

\blu{The most restricting requirement we make is therefore assumption (A1), which excludes interactions that do not have a hard-core component; notice how the condition $\big(\abs{x-y}\geq \HCdiam\big)$ can be generalised by $\big(y\notin x+\HardOpen\big)$, where $\HardOpen$ is any neighbourhood of the origin, for example an ellipse with minor axis equal to $\HCdiam$.}

The hard-core part of the interaction makes sure that any configuration has at most one point in each cube, if the cubes are small enough; this in turn allows us to derive the bound in \eqref{eq:th_uniqueness_region}.
\blu{Heuristically, the presence of the hard core means that the discretisation leads to a stationary approximating sequence. In the absence of it, the approximation would be only asymptotic.}
Many interactions, however, do not satisfy this assumption -- like the widely known Strauss pairwise interaction from \cite{strauss_1975}. In the case of a non-negative potential without hard core, as for the Strauss model with $\Potential (x,y) = \1_{\{\abs{x-y}\leq 1\}}$, we would still be able to apply the classical Dobrushin criterion and then consider the limit $\Width \to 0$. However, \blu{when taking the supremum in \eqref{eq:computation_dtv_removing_boundary} over all possible configurations $\BoundaryConf,\tilde \BoundaryConf$,} as in this case there is no restriction on the number of points in each cube, the uniqueness region we would obtain is $\Intensity <v_{\Dim}^{-1}$, where $v_{\Dim}$ is the volume of the unit ball.
Notice how this criterion does not depend on $\beta$ anymore.

One possible way to overcome this issue is to use the so-called Dobrushin--Pechersky uniqueness criterion \cite{dobrushin_pecherski_1983}. \blu{While the assumptions needed for this criterion are well understood (see \cite{pechersky_zhukov_1999}), the resulting conditions are too technical to obtain explicit values of the parameters.}

\subsubsection*{\blusec{Possible generalisations}}

It seems possible to develop new techniques inspired by the proof of Theorem \ref{thm:uniqueness_ours}. Heuristically, when taking $\Width \to 0$, it is less and less likely for a boundary condition to have more than one point in a small cube of side length $a$. One would then be able to make without the hard-core assumption by controlling instead the number of points in small boxes. This is the case for example with superstable potentials (\cite{ruelle_1970}) or with orderly point processes (\cite{daley_1974,ellis_1991}); in both cases one is able to bound the probability of having more than $2$ points in a given box in terms of its volume (\cite{schuhmacher_2005}).

By retracing the proof of Theorem \ref{thm:uniqueness_ours} with $\Ao\defeq\{\gamma:\#\gamma_{\Lambda_a}< 2\}$ in place of the $R$-hard-core configurations $\A$, we obtain the following bound, similar to \eqref{eq:tvestimate}:
	\begin{equation*}
    \forall \BoundaryConf, \tilde{\BoundaryConf}\in \Ao: \BoundaryConf_{\SuperLambda^c} = \tilde{\BoundaryConf}_{\SuperLambda^c}, \ d_{TV} \left( \Specification^{\Parameters}_{\SuperLambda, \BoundaryConf} , \Specification^{\Parameters}_{\SuperLambda, \tilde{\BoundaryConf}} \right) \leq \Intensity  \sup_{y \in \SuperLambda } \int_{\SuperLambda} \left( 1 - e^{- \beta \Potential (x,y) } \right) dx,
\end{equation*}
and we can set
\begin{equation*}
    k^{(a)'}_{i,j} \defeq \sup_{\substack{\BoundaryConf,\tilde\BoundaryConf\in\Ao\\\BoundaryConf_{\Lambda_{\Width,j}^c} = \tilde{\BoundaryConf}_{\Lambda_{\Width,j}^c}}} d_{TV} \left( \Specification^{\Parameters}_{\Lambda_{\Width, i}, \BoundaryConf}; \Specification^{\Parameters}_{\Lambda_{\Width, i}, \tilde{\BoundaryConf}} \right).
\end{equation*}
The point that remains open is whether one can apply the lattice Dobrushin criterion here. It is in fact not enough that $P(\lim_{a\to 0} \mathcal{A}^a)=1$, as what allows us to conclude the proof of Theorem \ref{thm:uniqueness_ours} is the existence of a \emph{positive} $a$, with $P(\mathcal{A}^a)=1$, that we can apply the classical Dobrushin criterion with (via Lemma \ref{lemma:dobrushin}). In this limiting case, the uniqueness of the Gibbs point process remains therefore a conjecture, which requires a not straightforward generalisation of Lemma \ref{lemma:dobrushin}.

\blu{Independently, M. Michelen and W. Perkins developed, in a recent preprint \cite{michelen-perkins}, a novel approach to prove uniqueness for repulsive and regular pair potentials, proposing a new explicit bound for the uniqueness region, by adapting the correlation decay method from theoretical computer science and using a recursive integral representation of the density of a point process.}

\subsection{Optimality of Theorem \ref{thm:uniqueness_ours} within the Dobrushin criterion}\label{sec:compareDob}

We believe the simplicity of the uniqueness region obtained in Theorem \ref{thm:uniqueness_ours}, as well as the ability to compare it to other criteria, is itself a justification for considering this method.
It is of course natural to ask whether taking the limit for $\Width\to 0$ actually yields a larger uniqueness region than what is obtained by just considering a fixed $\Width>0$, \blusec{i.e. when the second equality of assumption (A3) holds.}

\begin{figure}[ht]
\begin{minipage}[b]{.35\textwidth}
	\caption{Dobrushin uniqueness bound $a \mapsto \bar{\Intensity}_\beta(\Width)$ for the Hard-sphere model, for several values of $\Width\in[0,0.7]$. We find $\bar{\Intensity}_\beta(0)\simeq 0.318$.}\label{fig:dobrushin_uniqueness_hs}	
\end{minipage}\hfill
\begin{minipage}[b]{.6\textwidth}
	\centering
	\includegraphics[width=.8\textwidth]{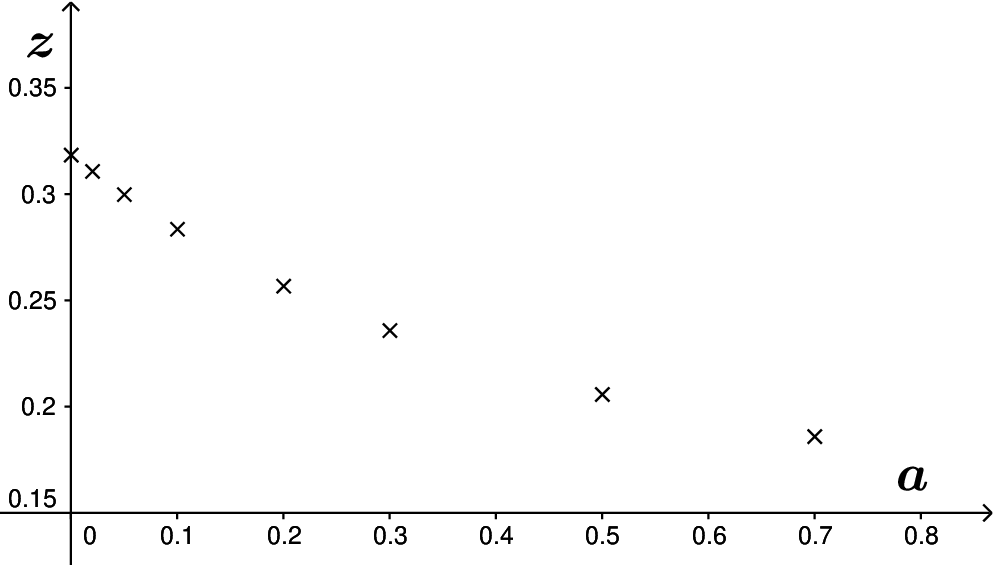}
\end{minipage}
\end{figure}

This translates to comparing the uniqueness regions obtained by fixing some $\Width>0$ and by taking the limit $a \to 0$, i.e. respectively (note that for $\bar\Intensity_\beta(a)$ we are using the possibly sub-optimal estimate of \eqref{eq:finalbound} for the total variation distance)
\blu{\begin{equation*}
	\bar\Intensity_\beta(\Width) \defeq \sup\,\{ \Intensity > 0:
	\Intensity \sup_{x \in \R^\Dim} \int_{\R^\Dim } \Psi_\Width(x,y)\, dy
	<1\} = \inf_{x\in\R^\Dim} \left(\int_{\R^\Dim } \Psi_\Width(x,y)\, dy\right)^{-1},\ a>0,
\end{equation*}
and
\begin{equation*}
	\bar\Intensity_\beta(0) \defeq \sup\, \{ \Intensity > 0: \Intensity\, \CRuelle < 1 \} = \left(\CRuelle\right)^{-1}.
\end{equation*}}
\blu{We note that, due to the form of the function $\Psi_a$ (defined in \eqref{eq:Psi}), the uniqueness region obtained for any fixed $a>0$ is of implicit nature -- while $\Psi_a$ is monotone in $a$, this monotonicity is broken when integrating and taking the supremum over all $x\in\R^d$ -- and it is therefore difficult to compare the uniqueness regions for different values of $a$.} We have performed, however, a numerical study in the case of the hard-sphere model $\Potential(x,y) = (+\infty) \1_{\{\abs{x-y}\leq 1\}}$ in dimension $\Dim=2$; Figure \ref{fig:dobrushin_uniqueness_hs} displays the upper bound $\bar\Intensity_{\beta}(\Width)$ of the Dobrushin uniqueness interval that we have computed numerically for several values of $a$, showing that taking the limit $\Width \to 0$ yields a larger Dobrushin uniqueness region than that of any fixed $\Width>0$. This indicates that the second equality of assumption (A3) holds.

\subsection{Comparison with the other uniqueness methods}
\label{sec:comparison}

In this subsection we briefly describe the two methods of cluster expansion and disagreement percolation, and see how the uniqueness regions that they yield compare to our Dobrushin criterion result. In particular, we provide a visual comparison between the three uniqueness methods, displayed below in Figure \ref{fig:comparison_methods_1} for two potentials with hard core, and in Figure \ref{fig:comparison_methods_strauss} for the case of the Strauss potential if the result of Theorem \ref{thm:uniqueness_ours} \blusec{holds without the hard core assumption.}

\begin{figure}[ht]
\begin{minipage}{\textwidth}
\begin{center}
\includegraphics[height=4.7cm]{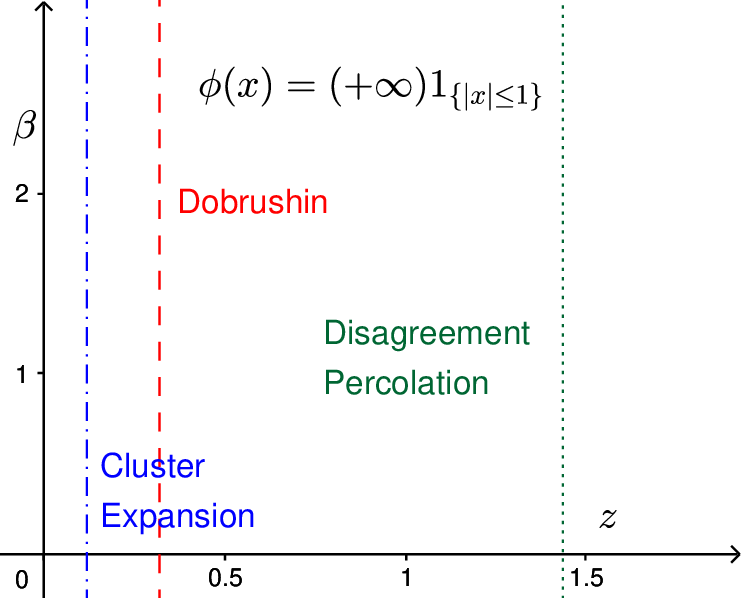}
\hspace{1cm}
\includegraphics[height=4.7cm]{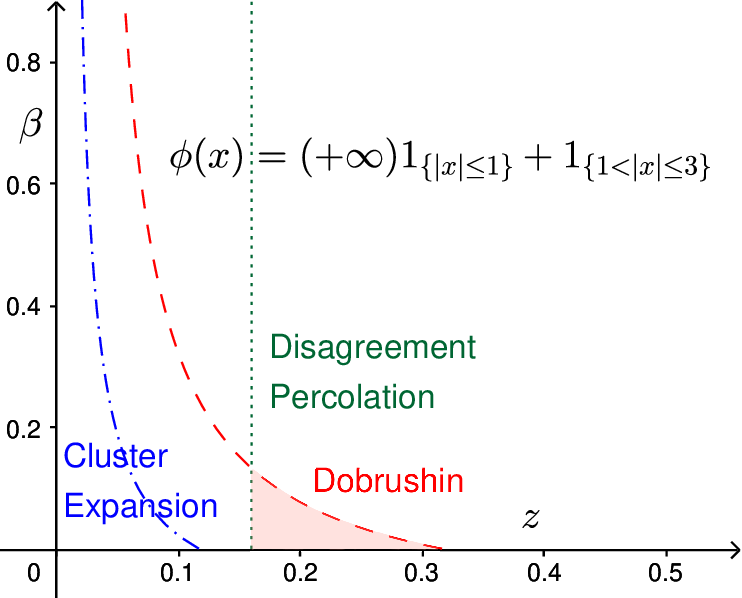}
\end{center}
\vspace{-1em}
\caption{The uniqueness regions (shaded where our method yields a larger bound) for the three different methods consist in the regions below the respective curves. Depicted here for the pure hard-core potential 
$\Potential(x,y)=(+\infty) \1_{\{\abs{x-y}\leq 1\}}$ on the left and for 
$\Potential(x,y)=(+\infty) \1_{\{\abs{x-y}\leq 1\}} + \1_{\{1<\abs{x-y}\leq 3\}}$ on the right.}
\label{fig:comparison_methods_1}
\vspace{.5em}
\end{minipage}
\begin{minipage}[b]{.5\textwidth}
	\caption{Uniqueness regions (shaded where our method yields a larger bound) for the Strauss potential $\Potential (x,y) = \1_{\{\abs{x-y}\leq 1\}}$, \blusec{assuming the result can be applied here.}}
\label{fig:comparison_methods_strauss}
\end{minipage}\hfill
\begin{minipage}[b]{.5\textwidth}
	\hspace{1.8em}\includegraphics[height=4.7cm]{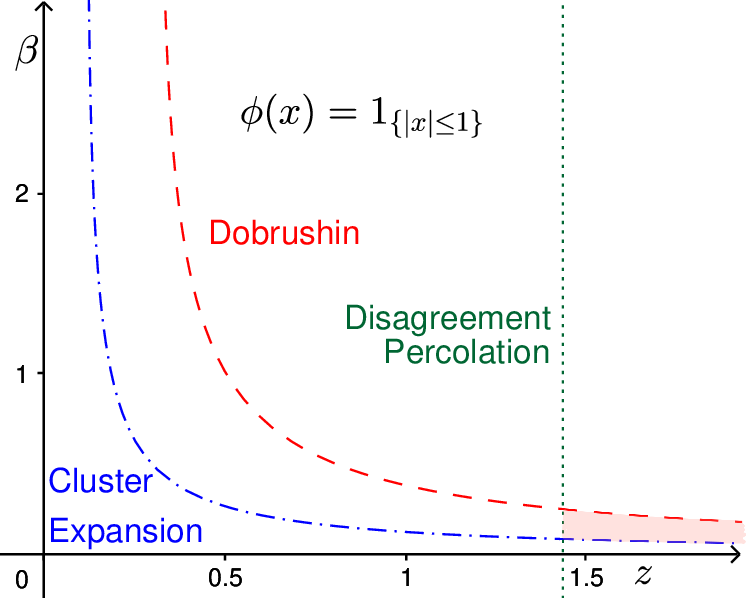}
\end{minipage}
\end{figure}

\subsubsection{Cluster expansion}\label{sec:cluster}

The method of cluster expansion was first developed for lattice systems in the 1980s (see e.g. \cite{Malyshev_1980}) and later extended to the continuous case. Indeed, different approaches exist to show convergence of the cluster terms (see e.g. \cite{malyshev_minlos_1991,poghosyan_ueltschi_2009,nehring_poghosyan_zessin_2012}).

For continuous point processes, the technique (from D. Ruelle, \cite{ruelle_livre_1969}) relies on a \emph{series expansion of the correlation functions}. More precisely, it consists in showing that the correlation functions of a Gibbs point process can be expressed as an absolutely converging series of cluster terms, and uniqueness is then proved by considering a set of integral equations -- the so-called \emph{Kirkwood--Salsburg equations} -- satisfied by the correlation functions, which can be reformulated as a fixed-point problem in an appropriately chosen Banach space, having therefore have a unique solution.

In \cite{jansen_2019}, S. Jansen presents a cluster expansion criterion for Gibbs point processes with a repulsive interaction $\phi$, \blu{with a condition similar to that presented by W.~G. Faris in \cite{Faris_2008} as a sufficient condition for the convergence of cluster expansion.} However -- as the most general form of this criterion yields an implicit uniqueness region -- in order to compare it with the domains obtained through the other methods, we only recall a specific uniqueness case that has an explicit region.

\begin{theorem}[\cite{jansen_2019}]
\label{thm:unique_cluster}
	Suppose there exist a non-negative measurable function $\mathsfbf{a}$ and some $t>0$ such that, for a.e. $x_0\in\R^\Dim$,
	\begin{equation*}
		z e^t\int \left(1- e^{-\beta\Potential(x_0,y)}\right)e^{\mathsfbf{a}(y)}dy\leq \mathsfbf{a}(x_0).
	\end{equation*}
	Then there exists a unique Gibbs measure $P\in\Gibbs$.
\end{theorem}
Restricting $\mathsfbf{a}$ to be a constant function, and remarking that $\max_{\mathsfbf{a}\geq 0} \mathsfbf{a}e^{-\mathsfbf{a}} = 1/e$, the condition of the above theorem holds for some $\mathsfbf{a}\geq 0$ and $t>0$ as soon as the classical Ruelle condition (e.g. Theorem 4.2.3 of \cite{ruelle_livre_1969}) holds
\begin{equation}
\label{eq:bound_cluster_expansion}
	z < \frac{1}{e}\, \CRuelle^{-1}.
\end{equation}
By considering a non-constant function $\mathsfbf{a}$ in Theorem \ref{thm:unique_cluster}, one could hope to obtain a better bound than \eqref{eq:bound_cluster_expansion}.

We also note that Fern\'andez, Procacci, and Scoppola were able to improve the classical cluster expansion bound using tree-graph estimates. Their approach was first presented in \cite{Fernandez_Procacci_2007}, where they also provide a comparison to the other methods; in \cite{Fernandez_Procacci_Scoppola_2007} the authors study the specific case of the $2$-dimensional hard-sphere model (where $\CRuelle = v_2$ is given by the volume of the unit ball in $\R^2$) and obtain from cluster expansion an improved bound which they estimate numerically as
\begin{equation*}
	z < 0.5107\, v_2^{-1}.
\end{equation*}
While larger than bound \eqref{eq:bound_cluster_expansion}, $z<e^{-1}v_2^{-1}$, this uniqueness region is still smaller than the one we obtain, i.e. $z<v_2^{-1}$.

\begin{remark}
\blu{The restriction to non-negative potentials greatly simplifies the cluster expansion approach. In this setting it is immediate to see that the correlation functions of any Gibbs point process are in the same Banach space where uniqueness holds: this is generally obtained by estimating the correlation functions via a so-called Ruelle bound, cf. \cite{ruelle_1971}, which is trivial in the case of repulsive interactions.}
\end{remark}

\subsubsection{Disagreement percolation}\label{sec:disperc}

The method of disagreement percolation was introduced for lattice systems by van den Berg and Maes \cite{vandenberg_1993,vandenberg_maes}.
The idea behind disagreement percolation is the construction of a coupling, sometimes called \emph{disagreement coupling}, which compares Gibbs specifications with two different boundary conditions outside of a given box, in such a way that the \emph{disagreement points} between the two Gibbs specification are ``connected'' to the boundary of the box.
If the probability of being connected to the boundary of an increasingly large box goes to zero, uniqueness holds.

Adapting the uniqueness result of \cite{Hofer-temmel_Houdebert_2018} to our setting reads:
\begin{theorem}[\cite{Hofer-temmel_Houdebert_2018}]
\label{thm:unique_perco}
Let $\Potential$ be a finite-range pair potential, i.e. there exists $r>0$ such that $\Potential(x,y)=0$ if $\abs{x-y}>r$. Then there exists at most one Gibbs point process $P\in\Gibbs$, for any $\beta\geq 0$ and for any activity
\begin{equation*}
	\Intensity < \frac{\Intensity_c(\Dim)}{r^\Dim},
\end{equation*}
where $\Intensity_c(\Dim)$ is the \emph{percolation threshold} of the Poisson--Boolean model in dimension $\Dim$ connecting points at distance at most one.
\end{theorem}
\begin{remark}
We recall that the \emph{Poisson--Boolean model} with fixed radius $1/2$ is a Poisson point process of balls $\pi^{z}$ where the configurations are collections of balls centred in the Poisson point in $\R^\Dim$ and with radius $1/2$. The percolation threshold $z_c(\Dim)$ is the smallest intensity value at which there exists $\pi^{z}$-almost surely an infinite connected component appears (i.e. the model \emph{percolates}).

Note that $\Intensity_c(\Dim)\geq v_\Dim^{-1}$ and, in the asymptotics as $\Dim\rightarrow\infty$, $\Intensity_c(\Dim)\sim v_\Dim^{-1}$, see \cite{penrose1996}.
We also remark that in  dimension $\Dim=2$ one can perform a simulation to see that $\Intensity_c(2) \simeq 1.43629$, see \cite{Quintanilla__Torquato_2000}.
\end{remark}
In \cite{Hofer-temmel_Houdebert_2018} the above bound is actually proved for interactions that do not necessarily come from a pair potential, without a hard core, and with a \emph{finite random range} that may depend on the (unbounded) marks of each point of the configuration.
However, it does not apply for infinite-range pair potentials, which are instead allowed in Theorem \ref{thm:uniqueness_ours}.

Finally, we remark that the disagreement percolation uniqueness region is independent of the parameter $\beta$ and depends only on the range of the interaction.
Typically, when $\beta$ is large, the disagreement percolation uniqueness region is larger than the one obtained from our bound in Theorem \ref{thm:uniqueness_ours}, but when $\beta$ is small, our uniqueness region is larger than the one coming from disagreement percolation (see the second potential of Figure \ref{fig:comparison_methods_1}).

\section*{Acknowledgments}
\noindent The authors would like to thank Sylvie R\oe lly: the final presentation of this work owes much to her.
We also thank the two referees for the many precise comments and remarks, which notably improved both the presentation and the content of the paper.

The research of the authors has been partially funded by the Deutsche Forschungsgemeinschaft (DFG)- Project-ID 318763901 - SFB1294.

\bibliographystyle{abbrv}

\end{document}